\documentclass[12pt]{amsart}
\usepackage{amsmath,amssymb,txfonts}
\usepackage{times}
\usepackage{bbm}
\usepackage[dvips]{graphicx}
\allowdisplaybreaks
\usepackage{amsxtra}
\usepackage[cp1252]{inputenc}

\usepackage{mathrsfs}

\numberwithin{equation}{section}

\newtheorem{theorem}{Theorem}[section]

\newtheorem{remark}[theorem]{Remark}

\newtheorem{lemma}[theorem]{Lemma}

\newcommand{\ct}[1]{\langle {#1}\rangle \lower.3ex\mathrm{$_{t}$}}
\newcommand{\lt}[1]{[ {#1}] \lower.3ex\mathrm{$_{t}$}}

\def\dint{\displaystyle\int}
\def\dsum{\displaystyle\sum}
\def\dsup{\displaystyle\sup}

\numberwithin{equation}{section}

\begin{document}
	
	\title{The $L^2$-boundedness of the variational Calder\'on-Zygmund operators}

\date{}

\maketitle

 \begin{center}
{\bf Yanping Chen \footnote{\small {Corresponding author.\ }}}\\
Department of Applied Mathematics,  School of Mathematics and Physics,\\
 University of Science and Technology Beijing,\\
Beijing 100083,  China \\
 E-mail: {\it yanpingch@ustb.edu.cn}
\vskip 0.1cm

{\bf Guixiang Hong}\\
School of Mathematics and Statistics\\
Wuhan University\\
Wuhan 430072, China\\
E-mail: {\it guixiang.hong@whu.edu.cn} \vskip 0.1cm

\end{center}

\renewcommand{\thefootnote}{}



\footnote{The project was in part supported by: Yanping Chen's
		National Natural Science Foundation of China (\# 11871096, \# 11471033); Guixiang Hong's National Natural Science Foundation of China ( \# 11601396).}

\renewcommand{\thefootnote}{\arabic{footnote}}
\setcounter{footnote}{0}

	\date{\today}

	\begin{abstract}
		In this paper, we verify the $L^2$-boundedness for  the jump functions and variations of Calder\'on-Zygmund singular integral operators with the underlying kernels satisfying
		\begin{align*}\int_{\varepsilon\leq |x-y|\leq N} K(x,y)dy=\int_{\varepsilon\leq |x-y|\leq N}K(x,y)dx=0\; \forall 0<\varepsilon\leq N<\infty,\end{align*}
		in addition to some proper size and smooth conditions. This result should be the first general criteria for the variational inequalities for kernels not necessarily of convolution type. The $L^2$-boundedness assumption that we verified here is also the starting point of the related results on the (sharp) weighted norm inequalities appeared in many recent papers.
	\end{abstract}
	
	\maketitle

\arraycolsep=1pt

 \section{Introduction}\label{s1}
\setcounter{equation}{0}

A singular integral operator in ${\Bbb R}^n$ is a continuous linear mapping $T$ from test functions
$\mathcal{D}({\Bbb R}^n)$ into distributions $\mathcal{D}'({\Bbb R}^n)$ associated to a Calder\'{o}n-
Zygmund standard kernel $K(x, y)$ in the sense that
\begin{align} Tf(x) =\int_{{\Bbb R}^n} K(x, y) f(y) \,dydx\end{align}whenever $f\in \mathcal{D}({\Bbb R}^n)$ and $x$ not in the support of $f$. If $T$ admits further a bounded extension on $L^2({\Bbb R}^n)$, that is $\|Tf\|_{L^2}\lesssim \|f\|_{L^2}, \;\forall f \in C_0^\infty({\Bbb R}^n),$ then $T$ is called a Calder\'{o}n-Zygmund operator associated with
the standard kernel $K.$
The $L^2$ boundedness criteria for singular integrals, commonly
known as $T1$ or $Tb$ theorems, which has arisen from  harmonic analysis and partial differential equations.
For $\varepsilon>0$, we define the truncated Calder\'on-Zygmund operator
$$T_\varepsilon f(x):=\int_{|x-y|>\varepsilon}K(x,y)f(y)dy.$$ The family of truncated Calder\'on-Zygmund operators $\{T_\varepsilon\}_{\varepsilon>0}$ will be denoted by $\mathcal {T}$ for simplicity. Let $2<q\leq\infty$, the $q$-variation operator for this family $\mathcal T=\{T_\varepsilon\}$ is defined by
\begin{align}V_q(\mathcal {T}f)(x):= V_q({T_{\varepsilon}}  f(x): \; \varepsilon>0 )=\sup_{\varepsilon_j\searrow 0}\bigg(\dsum_{j=1}^\infty|T_{\varepsilon_{j+1}}f(x)-T_{\varepsilon_{j}}f(x)|^q\bigg)^{1/q},\end{align}where the supremum is taken over all sequences $\{\varepsilon_j\}$ decreasing to zero. Note that when $q=\infty$, this is just the maximal Calder\'on-Zygmund operator.

Motivated by its studies in probability theory and ergodic theory \cite{Lep76, PiXu88, JKRW98, Bou89, AJS98, Gar, Qia98, JKRW98, JRW03, JRW00}, the first variational inequalities associated with singular integrals were established in \cite{CJRW2000} in the case of the Hilbert transform. Later, this result was generalized to higher dimensions for homogeneous singular integrals \cite{CJRW2002, JSW08, MaTo, CDHL}. The results for singular integrals with kernels of convolution type but without homogeneous properties seemingly first appeared in \cite{MST15}  under the additional assumption
$$\int_{\varepsilon\leq|x|\leq N}K(x)dx=0\;\forall 0<\varepsilon\leq N<\infty.$$

Even though this subject on variational inequalities has attracted a lot of attention from analysts and there appeared many papers \cite{GiTo04, LeXu2, DMT12, JSW08, MaTo, OSTTW12,  MiTr14, MST15, MTZ14, Zor14, KZ} ranging from the strengthened versions of Bourgain's variational estimates \cite{MTZ14} to the variational Carleson theorem \cite{OSTTW12} and the dimension-free variational estimates \cite{BMSW18}, there are very few \emph{complete} results on the Calder\'on-Zygmund operators associated with standard kernels of non-convolution type except the individual works like the one on Calder\'on's commutators \cite{CDHX} or the one on Cauchy integrals \cite{MaTo}. Indeed, in these papers  \cite{HLP, deZ16, DLY19, MTX1, MTX},  the $L^2$-variational inequality
 \begin{align}\label{assumption}
 \|V_q(\mathcal {T}f)\|_{L^2}\lesssim \|f\|_{L^2},\,\forall f\in L^2(\mathbb R^n),
 \end{align}
 has a priori  been assumed to complete their results  on the (sharp) weighted estimates of variational Calder\'on-Zygmund operators. Note that this assumption is quite natural since the $L^2$-boundedness of $\infty$-variational (or the maximal) Calder\'on-Zygmund operator are known to be true (see e.g. \cite{GR1, S1993})
 \begin{align}\label{maximal}
\|\sup_{\varepsilon>0}|T_\varepsilon f|\|_{L^2}\lesssim \|f\|_2.
\end{align}

However with a moment of thought, the assumption \eqref{assumption} for $q<\infty$ is wrong in the general setting because of the following two reasons. First, it is well-known that the validity of the $q$-variational inequalities with $q<\infty$ would imply immediately the a.e. convergence of truncated Calder\'on-Zygmund operators without the density argument; Secondly, there exist numerous examples such as the fractional singular integrals with complex-valued powers, which do not admit the a.e. convergence as $\varepsilon\rightarrow0$.

On the other hand, with the maximal inequality \eqref{maximal} at hand, by some compactness and density arguments,  there exists a decreasing subsequence $\{\varepsilon_j\}_j$ with $\varepsilon_j\rightarrow 0$ as $j\rightarrow\infty$ such that
\begin{align}\label{ass} \forall f\in L^2{(\mathbb R^n)},\ T_{\varepsilon_j}f\  \mathrm{converges}\  {a.e.} \; \mathrm{as} \ {j\rightarrow\infty}.\end{align}
One can find all the above assertions for instance in \cite{GR1, S1993}.

Thus the following conjecture for general Calder\'on-Zygmund operators arises naturally.

\smallskip

\noindent{\bf Conjecture.} Let $T$ be a Calder\'on-Zygmund operator with a standard kernel satisfying \eqref{ass} for a fixed sequence $\{\varepsilon_j\}_j$. Then for $2<q<\infty$, we have
\begin{align*}
\big\|V_q({T_{\varepsilon_j}}  f: \; j\in \mathbb N )\big\|_{L^2}
\lesssim \|f\|_{L^2}\;\forall f\in L^2{(\mathbb R^n)}.
\end{align*}

In the present paper, we give  a positive  answer to the conjecture for kernels satisfying
		$$\int_{\varepsilon\leq |x-y|\leq N}K(x,y)dy=\int_{\varepsilon\leq |x-y|\leq N}K(x,y)dx=0\; \forall 0<\varepsilon\leq N<\infty.$$ and verify this $L^2$-boundedness assumption \eqref{assumption} for a large class of Calder\'on-Zygmund operators with kernels not necessarily of convolution type.
		 Actually, we will prove a slightly stronger estimate, that is, the jump estimate. Let $\mathcal F=\{F_t:t\in\mathbb{R}_+\}$ be a family of Lebesgue measurable functions defined on $\mathbb{R}^n$.
Recall that   the $\lambda$-jump quantity $N_\lambda(\mathcal F)$ for $\lambda>0$ is defined as the supremum of $N$ over all increasing sequences $\{t_k\in \mathcal I\subset\mathbb{R}_+:0\leq k\leq N\}$ such that
$$ \min_{ k\in\{1,\dotsc, N\}}|F_{t_k}(x)-F_{t_{k-1}}(x)|>\lambda.$$
%Note that, the variation norm is just the maximum norm when $q=\infty$ and the finiteness of $q$-variation norm when $q<\infty$ implies the convergence along $\mathcal I$.

Recall that a standard Calder\'{o}n-Zygmund kernel  is a complex-valued function $K(x,y)$ defined on ${\Bbb R}^n\times{\Bbb R}^n\setminus\{(x,x): x\in {\Bbb R}^n\}$ satisfying
\begin{align} \label{k0}|K(x,y)|\lesssim \frac{1}{|x-y|^n}
\end{align}
 and there exists some $\theta\in (0,1)$ such that
for any $h\in {\Bbb R}^n$ with $2|h|\le |x-y|$,
\begin{align}\label{k2} |K(x+h,y)-K(x,y)|+|K(x,y+h)-K(x,y)|\lesssim \frac{|h|^\theta}{|x-y|^{n+\theta}}.
\end{align}

Motivated by the study of the sharp weighted normed inequalities in papers such as \cite{Ler17, deZ16, DLY19}, we will consider a weaker condition called the Dini condition
\begin{align}\label{k1} |K(x+h,y)-K(x,y)|+|K(x,y+h)-K(x,y)|\lesssim\frac{\omega(|h|/|x-y|)}{|x-y|^n},
\end{align}
where $\omega:[0,\infty)\rightarrow[0,\infty)$ with $\omega(0)=0$ is a function called a modulus of continuity of finite Dini norm
 \begin{align}\label{dini}
  \|\omega\|_{Dini}:=\int_0^1\omega(t)\frac{dt}{t}<\infty.
\end{align}
Let us recall some basic properties of a modulus of continuity which will be used implicitly in the paper. A modulus of continuity is sub-additive in the sense that
$$u\le t+s\Rightarrow \omega(u)\le \omega(t)+\omega(s).$$
Substituting $s=0$ one sees that $\omega(u)\le \omega(t)$ for all $0\le u\le t.$
It is easy to check that for any $c>0$ the integral in \eqref{dini} can be equivalently replaced by the sum over $2^{-j/c}$ with $j\in \Bbb N$ up to a $c$-dependent multiplicative constant.
The basic example is $\omega(t)=t^\theta$, which goes back to the Lipschitz condition with regularity index $\theta$.
Note that the composition or the sum of two modulus of continuity is again a modulus of continuity.
In particular, if $\omega(t)$ is a modulus of continuity and $\theta\in (0,1),$ then $\omega(t)^\theta$ and $\omega(t^\theta)$ are also moduli of continuity.

Now we present our main result.

\begin{theorem}\label{tJ2}
Let $K$ be a Calder\'on-Zygmund kernel satisfying \eqref{k0} and \eqref{k1} with some modulus of continuity $\omega$. If in addition, for any $0<\varepsilon\leq N$
\begin{equation}\label{tj2}\int_{\varepsilon \le |x-y|\leq N} K(x,y)\,dy=\int_{\varepsilon \le |x-y|\leq N} K(x,y)\,dx=0,\end{equation}
 then  \begin{align}\label{N12w}
\sup_{\lambda>0}\big\|\lambda\sqrt{N_\lambda(\mathcal{T}  f )}\big\|_{L^2}\lesssim (\|\omega^{1/2}\|_{Dini}^2+1)\|f\|_{L^2}\;\forall f\in L^2{(\mathbb R^n)}
\end{align}
and for all $2<q\leq\infty$,
\begin{align}\label{V}
\big\|V_q(\mathcal{T}  f )\big\|_{L^2}\lesssim (\|\omega^{1/2}\|_{Dini}^2+1)\|f\|_{L^2}\;\forall f\in L^2{(\mathbb R^n)}.
\end{align}
\end{theorem}
This result should be the first general criteria for the variational inequalities for kernels not necessarily of convolution type. As mentioned previously, the $L^2$-boundedness assumption that we verified here is also the starting point of the related results on the (sharp) weighted norm inequalities appeared in many recent papers, see e.g. \cite{HLP, deZ16, DLY19, MTX1, MTX} etc.

Regarding the proof of Theorem \ref{tJ2}, we will show only the jump estimate \eqref{N12w} since the variational estimate \eqref{V} can be proven similarly or can be deduced from the jump estimate by first obtaining all the $L^p$-variational estimates via Calder\'on-Zygmund theory as in \cite{MTX, MTX1} and then using the interpolation techniques as in Lemma 2.1 of \cite{JSW08}. Even though we shall show estimate \eqref{N12w} in Theorem \ref{tJ2} by checking separately the corresponding inequalities for the dyadic jump and the short variation as in most of the previously-cited papers (especially see Lemma 1.3 in \cite{JSW08}),  there appear a lot of difficulties. Indeed, for the dyadic jump estimate we have to establish Theorem \ref{tJ1} and two square function estimates \eqref{I2} and \eqref{I3}, all of which are results of $T1$ type;  while for the short variation the strategy in dealing with the kernels of convolution type fails completely in the present setting since the rapid decay estimates are not available, and our proof is based on again a $T1$ type argument.

The rest of this paper is organized as follows. In Section 2, we will show Lemma \ref{qj} and \ref{gj} which are key ingredients in dealing with the dyadic jump estimate. Then we establish Theorem \ref{tJ1} in Section 3. In Section 4 and Section 5,  we will show Theorem \ref{tJ2}. More precisely, in Section 4, we will show the following dyadic jump estimate,
\begin{align}\label{dyadic jump singular}
\sup_{\lambda>0}\big\|\lambda\sqrt{N_\lambda(\{ T_{2^k}f\}_{k\in \Bbb Z})}\big\|_{L^2}\lesssim\|f\|_{L^2};
\end{align}
while in Section 5, we will establish the short variational estimate,
\begin{align}\label{short variation singular}
\big\|S_2( \mathcal{T}f)\big\|_{ L^2}\lesssim\|f\|_{L^2}
\end{align}
with
$$
\begin{cases}
S_2( \mathcal{T}f)(x)=\bigg(\sum_{j\in\mathbb Z}[V_{2,j}(\mathcal{ T}f)(x)]^2\bigg)^{1/2};\\
V_{2,j}(\mathcal{ T}f)(x)=\bigg(\sup_
{2^j\leq t_0<\cdots<t_N<2^{j+1}}{\sum_{k=0}^{N-1}}| T_{t_{k+1}}f(x)- T_{t_k}f(x)|^2\bigg)^{1/2}.
\end{cases}
$$

\smallskip

\noindent{\bf Notation}. From now on, $p'=p/(p-1)$ represents the conjugate number of $p\in [1,\infty)$; $X\lesssim Y$ stands for $X\le C Y$ for a constant $C>0$ which is independent of the essential variables depending on $X\ \&\ Y$; and $X\simeq Y$ denotes  $X\lesssim Y\lesssim X$.

\section{Two key Lemmas}\label{s2}

Let $Q_s$ be defined by $Q_sf=\psi_s\ast f$ where $\psi_s(x)=s^{-n}\psi(x/s)$ with $\psi\in C_c^\infty(B(0,1))$ a radial real-valued function of mean zero. Noting that $Q_s$ is self-ajoint and \begin{align} \label{qs}\int_0^\infty Q_s^2\frac{ds}{s}=\mathcal{I},\end{align} where $\mathcal{I}$ is the identity operator on $L^2({\Bbb R}^n).$

\begin{lemma}\label{qj}
Let $\omega$ be a modulus of continuity and $T_j$ be a bounded linear operator in $L^2(\mathbb R^n)$ for each $j\in \Bbb Z$. Suppose
\begin{align}\label{t1}\|T_jQ_s\|_{L^2\rightarrow L^2}\lesssim \min\bigg\{\omega(\frac{2^j}{s}),\,\omega(\frac{s}{2^j})\bigg\}\end{align} and  \begin{align}\label{t2}\|Q_sT_j\|_{L^2\rightarrow L^2}\lesssim\min\bigg\{\omega(\frac{2^j}{s}),\,\omega(\frac{s}{2^j})\bigg\}.\end{align} Then we get
$$\bigg\|\sum_{j\in \Bbb Z}T_jf\bigg\|_{L^2}\lesssim \|\omega^{1/2}\|_{Dini}^2\|f\|_{L^2},\,\,\forall f\in L^2(\mathbb R^n),$$
where the infinite sum is convergent in $L^2$-norm.

\end{lemma}

{\emph{Proof}.}
There exists a function $g\in L^2$ with $\|g\|_{L^2}=1$ such that
\begin{align*}
\Big\|
\dsum_{j\in \Bbb Z}T_{j}f
\Big\|_{L^{2}}
&=\dsum_{j\in \Bbb Z}\int_{0}^{\infty}\int_{0}^{\infty}
\big\langle Q_{s}^{2}T_{j}Q_{t}^{2}f, g\big\rangle\frac{ds}{s}\frac{dt}{t}\\
&=\dsum_{j\in \Bbb Z}\int_{0}^{\infty}\int_{0}^{t}\big\langle Q^{2}_{s}T_{j}Q_{t}^{2}f, g\big\rangle
\frac{ds}{s}\frac{dt}{t}
+\dsum_{j\in \Bbb Z}\int_{0}^{\infty}
\int_{t}^{\infty}\big\langle Q_{s}^{2}T_{j}Q_{t}^{2}f, g\big\rangle
\frac{ds}{s}\frac{dt}{t}\\
&=:I+II.
\end{align*}
By symmetry, it is enough to consider the case $s\le t$,  so we only estimate $I.$
We split $I$ into three parts as follows:
\begin{align}\label{Tl}
I&=\dsum_{j\in \Bbb Z}\bigg(\int^{2^{j}}_{0}\int_{0}^{t}+\int_{2^{j}}^{\infty}\int_{0}^{2^{2j}t^{-1}}+\int_{2^{j}}^{\infty}\int^{t}_{2^{2j}t^{-1}}\bigg)\langle Q^{2}_{s}T_{j}Q_{t}^{2}f, g\rangle
\frac{ds}{s}\frac{dt}{t}\\&=:I_{1}+I_{2}+I_{3}\nonumber.
\end{align}
We now consider the first part appeared in the splitting of \eqref{Tl}. By the H\"{o}lder inequality, we get
 \begin{align*}
I_{1}&=\sum_{j\in \mathbb{Z}}\int_{0}^{2^{j}}\int_{0}^{t}\langle\omega^{-1/2}(s/2^j)Q_{s}T_{j}Q^{2}_{t}f, \omega^{1/2}(s/2^j)Q_{s}g\rangle\frac{ds}{s}\frac{dt}{t}\\
&\lesssim \bigg(\sum_{j\in \mathbb{Z}}\int_{0}^{2^{j}}\int_{0}^{t}\omega^{-1}(s/2^j)\|Q_{s}T_{j}Q^{2}_{t}f\|^{2}_{L^{2}} \frac{ds}{s}\frac{dt}{t}\bigg)^{\frac{1}{2}}\bigg(\sum_{j\in \mathbb{Z}}\int_{0}^{2^{j}}\int_{0}^{t}\omega(s/2^j)\|Q_{s}g\|^{2}_{L^{2}} \frac{ds}{s}\frac{dt}{t}\bigg)^{\frac{1}{2}}.
\end{align*}
Since $s\leq2^{j}$, by
$
\|Q_{s}T_{j}f\|_{L^{2}}\lesssim \omega(s/2^j)\|f\|_{L^{2}},
$ then \begin{align*}
I_{1}
&\lesssim \bigg(\int_{0}^{\infty}\sum_{2^{j}\geq t}\int_{0}^{t}\omega(s/2^j) \frac{ds}{s}\|Q^{2}_{t}f\|^{2}_{L^{2}}\frac{dt}{t}\bigg)^{\frac{1}{2}}\bigg(\int_{0}^{\infty}\int_{s}^{\infty}\sum_{2^{j}\geq t}\omega(s/2^j) \frac{dt}{t}\|Q_{s}g\|^{2}_{L^{2}}\frac{ds}{s}\bigg)^{\frac{1}{2}}.
\end{align*}
Since
\begin{align*}
\sum_{2^{j}\geq t}\int_{0}^{t}\omega(s/2^j) \frac{ds}{s}\le \Big(\sum_{2^{j}\geq t}\omega^{1/2}(t/2^j)\Big)\int_{0}^{1} \omega^{1/2}(s) \frac{ds}{s}\simeq \|\omega^{1/2}\|_{Dini}^2,
\end{align*} and
\begin{align*}
\int_{s}^{\infty}\sum_{2^{j}\geq t}\omega(s/2^j)\frac{dt}{t}\le  \int_{0}^{1}\Big(\sum_{2^{j}\geq t}\omega^{1/2}(t/2^j)\Big) \omega^{1/2}(t)\frac{dt}{t}\simeq \|\omega^{1/2}\|_{Dini}^2.
\end{align*}
Therefore,
\begin{align*}
I_{1}
&\lesssim  \|\omega^{1/2}\|_{Dini}^2\bigg(\int_{0}^{\infty}\|Q^{2}_{t}f\|^{2}_{L^{2}}\frac{dt}{t}\bigg)^{\frac{1}{2}}\|g\|_{L^{2}}\lesssim \|\omega^{1/2}\|_{Dini}^2\|f\|_{L^{2}}.
\end{align*}
Similarly, for $I_{2}$,  by $
\|Q_{s}T_{j}f\|_{L^{2}}\lesssim \omega(2^{-j}s)\|f\|_{L^{2}},$ we get \begin{align*}
I_{2}
&\lesssim\bigg(\int_{0}^{\infty}\int_{0}^{t} \sum_{2^{j}\geq s^{1/2}t^{1/2}} \omega(2^{-j}s)\|Q_{t}^{2}f\|^{2}_{L^{2}}\frac{ds}{s}\frac{dt}{t}\bigg)^{1/2}\bigg(\int_{0}^{\infty}\int_{s}^{\infty}\sum_{2^{j}\geq s^{1/2}t^{1/2}} \omega(2^{-j}s)\frac{dt}{t}\|Q_{s}g\|_{L^{2}}^2\frac{ds}{s}\bigg)^{1/2}.
\end{align*}
Note that $s\leq t$, $t\geq2^{j}$ combined with  $2^{2j}t^{-1}>s$ imply  $t\geq2^{j}\geq s^{1/2}t^{1/2}\geq s.$
Then \begin{align*}
\int_{s}^{\infty}\sum_{2^{j}\geq s^{1/2}t^{1/2}} \omega(2^{-j}s)\frac{dt}{t}&\le\Big(\sum_{2^{j}\geq s}\omega^{1/2}(2^{-j}s)\Big) \int_{s}^{\infty} \omega^{1/2}((\frac{s}{t})^{1/2})\frac{dt}{t}\simeq \|\omega^{1/2}\|_{Dini}^2,
\end{align*} and
\begin{align*}
\int_{0}^{t} \sum_{2^{j}\geq s^{1/2}t^{1/2}}\omega(2^{-j}s))\frac{ds}{s}&\le  \int_{0}^{t} \Big(\sum_{2^{j}\geq s}\omega^{1/2}(2^{-j}s)\Big)\omega^{1/2}((\frac{s}{t})^{1/2})\frac{ds}{s}\simeq \|\omega^{1/2}\|_{Dini}^2,
\end{align*}
then we get
\begin{align*}
I_{2}
&\lesssim\|\omega^{1/2}\|_{Dini}^2\bigg(\int_{0}^{\infty}\|Q_{t}f\|_{L^{2}}^{2}\frac{dt}{t}\bigg)^{1/2}
\bigg(\int_{0}^{\infty}\|Q_{s}g\|_{L^{2}}^2\frac{ds}{s}\bigg)^{1/2}\lesssim\|\omega^{1/2}\|_{Dini}^2\|f\|_{L^{2}}.
\end{align*}
Let us deal with $I_{3}$. Similarly, by  $\|T_jQ_tf\|_{L^2}\lesssim \omega({2^j}/{t})\|f\|_{L^2},$
we get
\begin{align*}
I_3&\lesssim\bigg(\int_{0}^{\infty}\int_{0}^{t}\sum_{2^{j}\leq s^{1/2}t^{1/2}}\omega({2^j}/{t})\|Q_{t}f\|_{L^{2}}^{2}\frac{ds}{s}\frac{dt}{t}\bigg)^{1/2}
\bigg( \int_{0}^{\infty}\int_{0}^{t} \sum_{2^{j}\leq s^{1/2}t^{1/2}}\omega({2^j}/{t})\|Q_{s}^2g\|^{2}_{L^{2}}\frac{ds}{s}\frac{dt}{t}\bigg)^{1/2}.
\end{align*}
Since \begin{align*}
\int_{0}^{t}\sum_{2^{j}\le s^{1/2}t^{1/2}}\omega({2^j}/{t})\frac{ds}{s}&\le \Big(\sum_{2^{j}\leq t}\omega^{1/2}(2^{j}/t)\Big) \int_{0}^{t} \omega^{1/2}((\frac{s}{t})^{1/2})\frac{ds}{s}\simeq \|\omega^{1/2}\|_{Dini}^2,
\end{align*} and
\begin{align*}
\int_{s}^{\infty} \sum_{2^{j}\le s^{1/2}t^{1/2}}\omega({2^j}/{t})\frac{dt}{t}&\le    \int_{s}^{\infty}\Big(\sum_{2^{j}\leq t}\omega^{1/2}(2^{j}/t)\Big) \omega^{1/2}((\frac{s}{t})^{1/2})\frac{dt}{t}\simeq \|\omega^{1/2}\|_{Dini}^2,
\end{align*}
then we get
\begin{align*}
I_3
\lesssim \|\omega^{1/2}\|_{Dini}^2\|f\|_{L^2}.
\end{align*}
Together with the estimates of $I_1,$ $I_2$ and $I_3$, we get
$$I\lesssim\|\omega^{1/2}\|_{Dini}^2\|f\|_{L^2}.$$\qed

\begin{remark}
Using the same proof, under the same assumptions as in Lemma \ref{qj}, one can show
$$\big\|\sum_{j\in \Bbb Z}\varepsilon_jT_jf\big\|_{L^2}\lesssim \|\omega^{1/2}\|_{Dini}^2\|f\|_{L^2}$$
uniformly for all independent Rademacher seqences $(\varepsilon_j)_j$. Then taking averages, we get the square function estimate
\begin{align}\label{s1}
\bigg\|\bigg(\sum_{j\in \Bbb Z}|T_jf|^2\bigg)^{1/2}\bigg\|_{L^2}\lesssim\|\omega^{1/2}\|^2_{Dini}\|f\|_{L^2}.
\end{align}
\end{remark}

However, we observe that the above square function requires less condition as in the following Lemma, which is crucial for the application in the proof of Theorem \ref{tJ2}.
\begin{lemma} \label{gj}
Let $\omega$ be a modulus of continuity and $T_j$ be a bounded linear operator in $L^2(\mathbb R^n)$ for each $j\in \Bbb Z$. If
\begin{align}\label{zj}
\|T_jQ_s\|_{L^{2}\rightarrow L^{2}}\lesssim\min\bigg\{\omega(\frac{2^j}{s}),\,\omega(\frac{s}{2^j})\bigg\},
\end{align} then we get
$$\bigg\|\bigg(\sum_{j\in \Bbb Z}|T_jf|^2\bigg)^{1/2}\bigg\|_{L^2}\lesssim\|\omega\|_{Dini}\|f\|_{L^2},\,\forall f\in L^2(\mathbb R^n).$$
\end{lemma}

\begin{remark}
Note that $\|\omega\|_{Dini}\lesssim \|\omega^{1/2}\|^2_{Dini}$, the square function estimate is better than \eqref{s1}.
\end{remark}

{\emph{Proof}.}  There exists $\{g_j\}\in L^2(\ell^2)$ with $\sum_{j\in \Bbb Z}\|g_j\|_{L^2}^2=1$ such that
\begin{align*}
\bigg\|\bigg(\sum_{j\in \Bbb Z}|T_{j}f|^{2}\bigg)^{\frac{1}{2}}\bigg\|_{L^{2}}
&=\sum_{j\in \Bbb Z}\int_{0}^{2^j}\langle T_{j}Q_{s}Q_{s}f,g_{j}\rangle\frac{ds}{s}
+\sum_{j\in \Bbb Z}\int_{2^j}^{\infty}\langle T_{j}Q_{s}Q_{s}f,g_{j}\rangle\frac{ds}{s}=:I+II.
\end{align*}
For $I$, using  the H\"{o}lder inequality  and $\|T_jQ_sf\|_{L^2}\lesssim\omega(\frac{s}{2^j})\|f\|_{L^2},$ we get
\begin{align*}
I&=\sum_{j\in \Bbb Z}\int_{0}^{2^j}\big\langle \omega^{-\frac{1}{2}}(\frac{s}{2^j})T_{j}Q_{s}Q_{s}f,\,\omega^{\frac{1}{2}}(\frac{s}{2^j})g_{j}\big\rangle\frac{ds}{s}\\
&\lesssim\bigg(\int_{0}^{\infty}\Big(\sum_{2^{j}>s}\omega(\frac{s}{2^j})\Big)\|Q_sf\|_{L^{2}}^{2}\frac{ds}{s}\bigg)^{1/2}
\bigg(\sum_{j\in \Bbb Z}\int_{0}^{2^j}\omega(\frac{s}{2^j})\|g_{j}\|_{L^{2}}^{2}\frac{ds}{s}\bigg)^{1/2}.
\end{align*}
Since
$
\sum_{2^{j}>s}\omega(\frac{s}{2^j})\simeq\|\omega\|_{Dini}
$
and $
\int_{0}^{2^j}\omega(\frac{s}{2^j})\frac{ds}{s}=\|\omega\|_{Dini},
$
then we get
\begin{align*}
I
&\lesssim\|\omega\|_{Dini}\bigg(\int_{0}^{\infty}\|Q_sf\|_{L^{2}}^{2}\frac{ds}{s}\bigg)^{1/2}
\bigg(\sum_{j\in \Bbb Z}\|g_{j}\|_{L^{2}}^{2}\bigg)^{1/2}\lesssim\|\omega\|_{Dini}\|f\|_{L^{2}}.
\end{align*}
For $II,$ similarly, using $\|T_jQ_sf\|_{L^2}\lesssim\omega(\frac{2^j}{s})\|f\|_{L^2},$ we get
\begin{align*}
II
&\lesssim\bigg(\int_{0}^{\infty}\sum_{2^j\leq s}\omega(\frac{2^j}{s})\|Q_sf\|_{L^{2}}^{2}\frac{ds}{s}\bigg)^{1/2}
\bigg(\sum_{j\in \Bbb Z}\int_{2^j}^{\infty}\omega(\frac{2^j}{s})\|g_{j}\|_{L^{2}}^{2}\frac{ds}{s}\bigg)^{1/2}\lesssim\|\omega\|_{Dini}\|f\|_{L^{2}}.
\end{align*}
Together with the estimates of $I$ and $II$, we finish the proof of Lemma \ref{gj}.\qed

\section{  A $T1$ type theorem }\label{s3}

\bigskip

In the course of establishing the main result, we show the following $T1$ type theorem, which might be of independent interest. Given a Calder\'on-Zygmund kernel $K$, we denote
$$\sigma_jf(x):=\int_{2^j\leq|x-y|< 2^{j+1}}K(x,y)f(y)dy,\;\forall j\in\mathbb Z.$$

\begin{theorem}\label{tJ1}
Let $K$ be a Calder\'on-Zygmund kernel satisfying \eqref{k0} and \eqref{k1} with some modulus of continuity $\omega$. If in addition, for any $j\in \Bbb Z,$
\begin{equation}\label{tj1}\int_{2^j\le |x-y|<2^{j+1}} K(x,y)\,dy=\int_{2^j\le |x-y|<2^{j+1}} K(x,y)\,dx=0,\end{equation}
 then $\sum^M_{j=-N}\sigma_j$ converges in the strong operator topology  as $M,N\rightarrow\infty$.
 Moreover the limit denoted by $T$ satisfies
\begin{equation}
\|Tf\|_{L^2}\lesssim (\|\omega^{1/2}\|_{Dini}^2+1)\|f\|_{L^2} \;\forall f\in L^2(\mathbb R^n).
\end{equation}
\end{theorem}
 When the Dini condition is replaced by the Lipschitz condition, that is, $\omega(t)=t^{\theta}$ with $\theta\in (0,1)$, Theorem \ref{tJ1} has been previously established, see e.g. Proposition 8.5.3 of \cite{GR1}.

{\emph{Proof}.} We split $K(x,y)$ into $K(x,y)=\sum_{j\in \Bbb Z}K_j(x,y),$ where $K_j(x,y)=K(x,y)\chi_{\{2^j\le |x-y|< 2^{j+1}\}}.$ Then $Tf=\sum_{j\in \Bbb Z}\sigma_jf,$ where $$\sigma_jf(x):=\int_{{\Bbb R}^n} K_j(x,yf(y)\,dy.$$
 By \eqref{k0}, it is easy to verify that for any fixed $j\in \Bbb Z,$\begin{align}\label{kj0}
 |K_{j}(x,y)|&\lesssim \frac{1}{|x-y|^n}\chi_{\{2^j\le|x-y|< 2^{j+1}\}},
\end{align}
and by \eqref{k1} for $2|h|\le 2^j,$ we get
\begin{align*}
&|K_{j}(x,y+h)-K_{j}(x,y)|\\
&\le|K(x,y+h)-K(x,y)|\chi_{\{2^{j}\le|x-y-h|<2^{j+1}\}}
+|K(x,y)||\chi_{\{2^{j}\le |x-y-h|<2^{j+1}\}}-\chi_{\{2^{j}\le|x-y|<2^{j+1}\}}|\\
&\lesssim\frac{\omega(|h|/|x-y|)}{|x-y|^n}\chi_{\{2^{j}/2\le |x-y|\le 2^{j+2}\}}
+\frac{1}{|x-y|^n}|\chi_{\{2^{j}\le|x-y-h|<2^{j+1}\}}-\chi_{\{2^{j}\le|x-y|<2^{j+1}\}}|.
\end{align*}
 Note that
\begin{align*}
\chi_{\{2^{j}\le|x-y-h|<2^{j+1}\}}-\chi_{\{2^{j}\le|x-y|<2^{j+1}\}}\neq0,
\end{align*}
if and only if at least one of the following two statements holds:

(i)\, $2^{j}\le|x-y-h|<2^{j+1}$ and $\{|x-y|<2^{j}\}\cup\{|x-y|\ge2^{j+1}\}$;

(ii)\, $2^{j}\le|x-y|<2^{j+1}$ and $\{|x-y-h|< 2^{j}\}\cup\{|x-y-h|\ge2^{j+1}\}$;

{\noindent This together with the fact that $2|h|\leq 2^j$ implies the following two cases:}

(i)\, $2^{j}-|h|\leq|x-y|\leq2^{j}+|h|;$

(ii)\,  $2^{j+1}-|h|\leq|x-y|\leq2^{j+1}+|h|;$

{\noindent Therefore by \eqref{k0} we get for $2|h|\le 2^j$ }
\begin{align}\label{kj1}
&|K_{j}(x,y+h)-K_{j}(x,y)|\\
&\lesssim\frac{\omega(|h|/|x-y|)}{|x-y|^n}\chi_{\{2^{j}/2\le |x-y|\le 2^{j+2}\}}
+\frac{1}{2^{jn}}\chi_{\{2^{j}-|h|\leq|x-y|\leq2^{j}+|h|\}\cup\{2^{j+1}-|h|\leq|x-y|\leq2^{j+1}+|h|\}}.\nonumber
\end{align}
Similarly, we get
 for $2|h|\le 2^j,$
\begin{align}\label{kj2}
&|K_{j}(x+h,y)-K_{j}(x,y)|\\
&\lesssim\frac{\omega(|h|/|x-y|)}{|x-y|^n}\chi_{\{2^{j}/2\le |x-y|\le 2^{j+2}\}}
+\frac{1}{2^{jn}}\chi_{\{2^{j}-|h|\leq|x-y|\leq2^{j}+|h|\}\cup\{2^{j+1}-|h|\leq|x-y|\leq2^{j+1}+|h|\}}.\nonumber
\end{align}

   In the following we would like to  use \eqref{kj0}, \eqref{kj1}  and $\sigma_j1=0$ (resp. \eqref{kj0}, \eqref{kj2}  and $\sigma_j^*1=0$) for any $j\in \Bbb Z$, to prove that
    \begin{align}\label{t3}\|\sigma_jQ_sf\|_{L^2}\lesssim \min\{\omega_1(\frac{2^j}{s}),\,\omega_1(\frac{s}{2^j})\}\|f\|_{L^2}, (\mathrm{resp.}\;\|Q_s\sigma_jf\|_{L^2}\lesssim \min\big\{\omega_1(\frac{2^j}{s}),\,\omega_1(\frac{s}{2^j})\big\}\|f\|_{L^2})
    \end{align}
     where $\omega_1(t)=\omega(t)+t^\theta$, with some $\theta\in (0,1)$ and
$ \|\omega_1^{1/2}\|_{Dini}\lesssim \|\omega^{1/2}\|_{Dini}+1.
$
Then by Lemma \ref{qj}, we get \begin{align}\label{t4}\|Tf\|_{L^2}\lesssim (\|\omega^{1/2}\|_{Dini}^2+1)\|f\|_{L^2}.\end{align}

 First, we are ready to verify the first estimate in \eqref{t3}.
First of all, for the case of $2^{j-1}\le s,$ by $\sigma_j1=0$ for any fixed $j\in \Bbb Z,$ and \eqref{kj0}
\begin{align*}
|\sigma_{j}Q_{s}f(x)|&=\bigg|\int_{{\Bbb R}^n}\int_{{\Bbb R}^n} K_{j}(x,z)\big(\psi_{s}(z-y)-\psi_{s}(x-y)\big)dzf(y)dy\bigg|\\&\lesssim\int_{|x-y|\leq4s}\int_{|x-z|\le s} \frac{2^{j/2}}{|x-z|^{n+1/2}}|\psi_{s}(z-y)-\psi_{s}(x-y)|dz|f(y)|dy\\&\quad+\int_{{\Bbb R}^n}\int_{|x-z|> s} \frac{2^{j/2}}{|x-z|^{n+1/2}}|\psi_{s}(z-y)|dz|f(y)|dy+\int_{{\Bbb R}^n}\int_{|x-z|> s}\frac{2^{j/2}}{|x-z|^{n+1/2}}|\psi_{s}(x-y)|dz|f(y)|dy\\
&\lesssim\int_{|x-y|\leq4s}\int_{|x-z|\le s}\frac{2^{j/2}}{|x-z|^{n+1/2}}\frac{|x-z|}{s^{n+1}}dz|f(y)|dy+(\frac{2^{j}}{s})^{1/2}(M^2f(x)+Mf(x))\\
&\lesssim(\frac{2^{j}}{s})^{1/2}(Mf(x)+M^2f(x)).
\end{align*} Here $Mf(x)$ is the usual Hardy-Littlewood maximal function and $M^2f=M(Mf)$.
Then we get
\begin{align}\label{ts}
\|\sigma_{j}Q_{s}f\|_{L^2}\lesssim (\frac{2^{j}}{s})^{1/2}\|f\|_{L^2}.
\end{align}
Next, for $s\leq2^{j-1},$
using  $Q_{s}1=0$ and the estimate \eqref{kj1}, we get
\begin{align*}
|\sigma_{j}Q_{s}f(x)|&=\bigg|\int_{{\Bbb R}^n}\int_{{\Bbb R}^n}(K_{j}(x,z)-K_{j}(x,y))\psi_{s}(z-y)\,dzf(y)\,dy\bigg|\\
&\lesssim\int_{{\Bbb R}^n}\int_{{\Bbb R}^n}\frac{\omega(|z-y|/|x-y|)}{|x-y|^n}\chi_{\{2^{j}/2\le |x-y|\le 2^{j+2}\}}|\psi_{s}(z-y)|\,dz|f(y)|\,dy\\
&\quad+\frac{1}{2^{jn}}\int_{\{2^j-s\le |x-y|\le 2^j+s\}\cup\{2^{j+1}-s\leq|x-y|\leq2^{j+1}+s\}}\int_{{\Bbb R}^n} |\psi_{s}(z-y)|\,dz|f(y)|\,dy\\&\lesssim  \omega(\frac{s}{2^j}) Mf(x)+ \frac{1}{2^{jn}}\bigg(\int_{\{2^j-s\le |x-y|\le 2^j+s\}\cup\{2^{j+1}-s\leq|x-y|\leq2^{j+1}+s\}}\,dy\bigg)^{1/p'}\bigg(\int_{|x-y|\le 2^{j+1}}|f(y)|^p\,dy\bigg)^{1/p}\\
&\lesssim \omega(\frac{s}{2^j}) Mf(x)+\frac{1}{2^{jn}}2^{\frac{j(n-1)}{p'}}s^{1/p'}\bigg(\int_{|x-y|\le 2^{j+1}}|f(y)|^p\,dy\bigg)^{1/p}\\
&\lesssim \omega(\frac{s}{2^j}) Mf(x)+(\frac{s}{2^j})^{1/p'}(M(|f|^p))^{1/p}(x),
\end{align*}
where $1<p<2$ can be chosen arbitrarily.
Then we get for some $\theta\in (0,1),$\begin{align*}
\|\sigma_{j}Q_{s}f\|_{L^{2}}&\lesssim\big(\omega(\frac{s}{2^j})+\big(\frac{s}{2^{j}}\big)^{\theta}\big)\|f\|_{L^{2}}.
\end{align*}
This together \eqref{ts} implies the desired estimate \eqref{t3}.

Now we are ready to prove the second estimate in \eqref{t3}.
First, for the case of $2^{j-1}\le s,$ by  $\sigma_j^*1=0$ for any fixed $j\in \Bbb Z,$  \eqref{kj0} and the smoothness of $\psi$, we get
\begin{align*}
|Q_{s}\sigma_jf(x)|&=\bigg|\int_{{\Bbb R}^n}\int_{{\Bbb R}^n} \big(\psi_{s}(x-z)-\psi_{s}(x-y)\big) K_{j}(z,y)dzf(y)dy\bigg|\\&\lesssim \int_{{\Bbb R}^n}\int_{|z-y|\le s} |\psi_{s}(x-z)-\psi_{s}(x-y)| |K_{j}(z,y)|dz|f(y)|dy\\&\quad+\int_{{\Bbb R}^n}\int_{|z-y|> s}\frac{2^j}{|z-y|^{n+1}}|f(y)|dy|\psi_{s}(x-z)|dz+\int_{{\Bbb R}^n}\int_{|z-y|> s}\frac{2^j}{|z-y|^{n+1}}dz|\psi_{s}(x-y)||f(y)|dy\\
&\lesssim(\frac{2^{j}}{s})^{1/2}\frac{1}{s^{n}}\int_{|x-y|\leq4s}\int_{|z-y|\le s}\frac{2^{j/2}|z-y|^{1/2}}{(2^j+|z-y|)^{n+1}}dz|f(y)|dy+\frac{2^{j}}{s}(M^2f(x)+Mf(x))\\
&\lesssim(\frac{2^{j}}{s})^{1/2}(Mf(x)+M^2f(x)).
\end{align*}
Then we get
\begin{align}\label{ts2}
\|Q_{s}\sigma_{j}f\|_{L^2}\lesssim (\frac{2^j}{s})^{1/2}\|f\|_{L^2}.
\end{align}
Thus it remains to show the case of  $s\le 2^{j-1}.$
Since $Q_{s}1=0$ for any fixed $s$, then by \eqref{kj2}, as the arguments after \eqref{ts}, one gets
\begin{align*}
|Q_{s}\sigma_{j}f(x)|
&=\bigg|\int_{{\Bbb R}^n}\int_{{\Bbb R}^n}\psi_{s}(x-z)\big(K_{j}(z,y)-K_{j}(x,y)\big)\,dzf(y)\,dy\bigg|\\
&\lesssim \int_{2^{j}/2\le |x-y|\le 2^{j+2}}\int_{{\Bbb R}^n}\frac{\omega(|x-z|/|x-y|)}{|x-y|^n}|\psi_{s}(x-z)|\,dz|f(y)|\,dy\\
&\quad+\frac{1}{2^{jn}}\int_{\{2^j-s\le |x-y|\le 2^j+s\}\cup\{2^{j+1}-s\le |x-y|\le 2^{j+1}+s\}}\int_{{\Bbb R}^n} |\psi_{s}(x-z)|\,dz|f(y)|\,dy\\
&\lesssim \omega(\frac{s}{2^j}) Mf(x)+(\frac{s}{2^j})^{1/p'}(M(|f|^p))^{1/p}(x),
\end{align*}
 where $1<p<2$ can be chosen arbitrarily.
Then we get for some $\theta\in (0,1),$\begin{align*}
\|Q_{s}\sigma_{j}f\|_{L^{2}}&\lesssim \big(\omega(\frac{s}{2^j})+\big(\frac{s}{2^{j}}\big)^{\theta})\|f\|_{L^{2}}.
\end{align*} This together with \eqref{ts2} implies the second estimate in \eqref{t3}.\qed

\section{Proof of  Theorem \ref{tJ2}---The dyadic jump estimate \eqref{dyadic jump singular}}\label{s5}

Now choose a radial $\varphi\in C^{\infty}_{c}(B(0,1/2))$ with $\int\varphi=1$, and let $\varphi_{j}(x)= 2^{-jn}\varphi(2^{-j}x)$ for $j\in \Bbb Z$. Following \cite{DR86}, we write
$T_{j}f(x)=\int_{|x-y|>2^{j+1}}K(x,y)f(y)\,dy$ and $T^{j}f(x)=\int_{|x-y|\leq 2^{j+1}}K(x,y)f(y)\,dy$, then
\begin{align*}
T_{j}f(x)
&=\varphi_{j}\ast Tf-\varphi_{j}\ast T^{j}f+(\delta-\varphi_{j})\ast T_{j}f=:I_1-I_2+I_3.
\end{align*}
In the following, we will estimate $I_1,$ $I_2$ and $I_3$, respectively.

\subsection{{Proof of $I_1$}} To estimate the term $I_1,$ we first give a known lemma which is stated as follows.
\begin{lemma}\label{Phi} (\cite{JSW08},\,\cite{DHL})  Let $\varphi\in C^\infty_c ({\Bbb R}^n)$ with $\int \varphi=1.$ Then for $1<p<\infty$,
$$
\sup_{\lambda>0}\Big\|\lambda\sqrt{ N_\lambda (\{\varphi_j\ast f\}_{j\in\mathbb Z})}\Big\|_{L^p}\lesssim\|f\|_{L^p} \;\forall f\in L^p(\mathbb R^n).
$$
\end{lemma}
For $I_{1}$,
 by Lemma \ref{Phi} and Theorem \ref{tJ1}, it is easy to get\begin{align*}
\sup_{\lambda>0}\Big\|\lambda \sqrt{N_{\lambda}(\{\varphi_{j}\ast Tf\}_{j\in \Bbb Z})}\Big\|_{L^2}\lesssim \|Tf\|_{L^{2}}\lesssim (\|\omega^{1/2}\|_{Dini}^2+1)\|f\|_{L^{2}}.
\end{align*}

\subsection{{Proof of $I_2$}} For $I_{2}$, by the definition of the jump quantity and Chebychev's inequality,
\begin{align*}
\sup_{\lambda>0}\Big\|\lambda \sqrt{N_{\lambda}(\{\varphi_{j}\ast T^{j}f)\}_{j\in \Bbb Z}}\Big\|_{L^2}\lesssim  \bigg\|\bigg(\sum_{j\in\mathbb{Z}}|\varphi_{j}\ast T^{j}f|^{2}\bigg)^{\frac{1}{2}}\bigg\|_{L^{2}}.
\end{align*}
So to estimate $I_2,$ we need to show  \begin{align} \label{I2}\bigg\|\Big(\sum_{j\in\mathbb{Z}}|\varphi_{j}\ast T^{j}f|^{2}\Big)^{\frac{1}{2}}\bigg\|_{L^{2}} \lesssim (\|\omega^{1/2}\|_{Dini}^2+1)\|f\|_{L^2}. \end{align} For any $j\in \Bbb Z,$ set $A_{j}:=\varphi_{j}\ast T^{j}$ and write the kernel of $A_{j}$ as $$a_{j}(x,y)=\int_{{\Bbb R}^n}\varphi_{j}(x-z)K^j(z,y)\,dz,$$ where $K^j(z,y)=K(z,y)\chi_{\{|z-y|\le 2^{j+1}\}}$. Observe that $a_{j}$ is supported on the set $\{|x-y|\lesssim2^{j}\}$.

We will apply Lemma \ref{gj} to handle the inequality \eqref{I2}. By Lemma \ref{gj}, we need to verify that $A_j$ satisfies \eqref{zj}. First, we deal with the case of $2^j\le s.$
Since $\int_{{\Bbb R}^n} K^{j}(z,y)\,dz=\lim_{\varepsilon\rightarrow0}\int_{|z-y|>\varepsilon} K^{j}(z,y)\,dz=0,$ we get
$
a_{j}(x,y)=\int_{{\Bbb R}^n}(\varphi_{j}(x-z)-\varphi_{j}(x-y))K^j(z,y)\,dz.
$
 A trivial computation gives that
\begin{align}\label{aj}
|a_{j}(x,y)|\lesssim\int_{|z-y|\leq 2^{j+1}}\frac{|y-z|}{2^{j(n+1)}}|K(z,y)|\,dz\lesssim\frac{1}{2^{jn}}\chi_{\{|x-y|\lesssim2^{j}\}}.
\end{align}
Since $\int_{{\Bbb R}^n} K^{j}(\eta,z)\,dz=\lim_{\varepsilon\rightarrow0}\int_{|\eta-z|>\varepsilon} K^{j}(\eta,z)\,dz=0,$ then $\int_{{\Bbb R}^n} a_j(x,z)\,dz=\int_{{\Bbb R}^n} \int_{{\Bbb R}^n}\varphi_{j}(x-\eta)K^j(\eta,z)\,d\eta\,dz=0,$ we get
\begin{align*}
 |A_{j}Q_{s}f(x)|&=\bigg|\int_{{\Bbb R}^n}\int_{{\Bbb R}^n} a_{j}(x,z)\big(\psi_{s}(z-y)-\psi_{s}(x-y)\big)dzf(y)dy\bigg|\\
 &\lesssim \int_{|x-y|\lesssim s}\int_{|x-z|\lesssim2^j} \frac{1}{2^{jn}} \frac{|z-x|}{s^{n+1}} dz|f(y)|dy\lesssim\frac{2^{j}}{s}Mf(x).
\end{align*}
Thus,
\begin{align}\label{A}
 \|A_{j}Q_{s}f\|_{L^2}\lesssim\frac{2^{j}}{s}\|f\|_{L^2}.
\end{align}
Next, we consider the case of $s\leq 2^j.$ For  $|h|\leq2^j$, since $\int_{{\Bbb R}^n} K^{j}(z,y)\,dz=0$, we get
\begin{align}\label{ajh}
&|a_j(x,y+h)-a_j(x,y)|\\\nonumber
&=\bigg|\int_{{\Bbb R}^n}\big(\varphi_j(x-z)-\varphi_j(x-y-h)\big)K^j(z,y+h)\,dz-
\int_{{\Bbb R}^n}\big(\varphi_j(x-z)-\varphi_j(x-y)\big)K^j(z,y)\,dz\bigg|.\nonumber
\end{align}
Deduced from the support of $\varphi_j$ and $K^j$,  one has  $a_j(\cdot,\cdot+h)-a_j(\cdot,\cdot)$ is supported on the subset $\{|x-y|\lesssim2^{j}\}.$ We proceed with the proof by separating the domain of integration in \eqref{ajh},
\begin{align*}
&|a_j(x,y+h)-a_j(x,y)|\\
&=\bigg|\int_{|z-y|<2|h|}\big(\varphi_j(x-z)-\varphi_j(x-y-h)\big)K^j(z,y+h)\,dz-
\int_{|z-y|<2|h|}\big(\varphi_j(x-z)-\varphi_j(x-y)\big)K^j(z,y)\,dz\bigg|\\&\quad+\bigg|\int_{|z-y|\ge 2|h|}\big(\varphi_j(x-z)-\varphi_j(x-y-h)\big)K^j(z,y+h)\,dz-
\int_{|z-y|\ge 2|h|}\big(\varphi_j(x-z)-\varphi_j(x-y)\big)K^j(z,y)\,dz\bigg|\\&=:A_1+A_2.
\end{align*}
For $A_1,$ by \eqref{k0}, we get $|K^j(z,y)|\lesssim \frac{1}{|z-y|^n}\chi_{\{|z-y|\le 2^{j+1}\}},$ then for some $\theta\in (0,1),$
\begin{align*}
A_1
&\lesssim\int_{\{|z-y-h|\le 2^{j+1}\}\cap \{|z-y|<2|h|\}}\frac{|z-y-h|}{2^{j(n+1)}}\frac{1}{|z-y-h|^n}\,dz+\int_{\{|z-y|\le 2^{j+1}\}\cap \{|z-y|<2|h|\}}\frac{|z-y|}{2^{j(n+1)}}\frac{1}{|z-y|^n}\,dz\\ &\lesssim\int_{|z-y-h|\le 2^{j+1}}\frac{|h|^\theta}{2^{j(n+1)}}\frac{1}{|z-y-h|^{n+\theta-1}}\,dz+\int_{||z-y|\le 2^{j+1}}\frac{|h|^\theta}{2^{j(n+1)}}\frac{1}{|z-y|^{n+\theta-1}}\,dz\\&\lesssim \frac{|h|^\theta}{2^{j(n+\theta)}}\chi_{\{|x-y|\lesssim2^j\}}.
\end{align*}
For $A_2,$ we get
\begin{align*}
A_2& \lesssim
\int_{|z-y|\ge 2|h|}|\varphi_j(x-y)-\varphi_j(x-y-h)||K^j(y,z)|\,dz\\&\quad+\int_{|z-y|\ge 2|h|}|\varphi_j(x-z)-\varphi_j(x-y-h)|
|[K(y,z+h)-K(y,z)]\chi_{\{|y-z-h|\le 2^{j+1}\}}|\,dz\\
&\quad+\int_{|z-y|\ge 2|h|}|\varphi_j(x-z)-\varphi_j(x-y-h)|
|K(y,z)\big(\chi_{\{|y-z-h|\le 2^{j+1}\}}-\chi_{\{|y-z|\le 2^{j+1}\}}\big)|\,dz\\
&=:A_{3}+A_{4}+A_5.
\end{align*}
For $A_{3}$, by $|K^j(z,y)|\lesssim \frac{1}{|z-y|^n}\chi_{\{|z-y|\le 2^{j+1}\}},$ we get for  $|h|\le 2^j$,
\begin{align*}
A_{3}
&\lesssim\int_{2|h|\le |y-z|\leq2^{j+1}}\frac{|h|}{2^{j(n+1)}}\cdot\frac{1}{|y-z|^n}\,dz\lesssim\frac{|h|^\theta}{2^{j(n+\theta)}}\chi_{\{|x-y|\lesssim2^j\}}.
\end{align*}
For $A_{4}$, since $|z-y|\geq2h$, one has $\frac{|y-z|}{2}<|y-z-h|\le 2|y-z|$.  By the sub-additivity of $\omega(t)$, we get for any fixed $k\in \Bbb N_+,$  $\omega(kt)\le k\omega(t).$ Then apply \eqref{k1},
\begin{align*}
A_{4}&\lesssim\int_{|h|\le|y-z|\le2^{j+2}}\frac{|z-y-h|}{2^{j(n+1)}}\cdot\frac{\omega(|h|/|y-z|)}{|y-z|^{n}}\,dz\\
&\lesssim\frac{\omega^{1/2}(1)}{2^{j(n+1)}}\sum_{k=-\infty}^1\omega^{1/2}(2^{-k}|h|/2^{j})\int_{2^{k+j}<|y-z|\le2^{k+j+1}}\frac{1}{|y-z|^{n-1}}\,dz\\
&\lesssim\frac{\omega^{1/2}(1)}{2^{j(n+1)}}\sum_{k=-\infty}^12^{k+j}2^{-k/2}\omega^{1/2}(|h|/2^{j})\\
&\lesssim \|\omega^{1/2}\|_{Dini}\frac{\omega^{1/2}(|h|/2^{j})}{2^{jn}}\chi_{\{|x-y|\lesssim2^j\}}.
\end{align*}
To estimate $A_5$, note that
$$\chi_{\{|y-z-h|\le 2^{j+1}\}}-\chi_{\{|y-z|\le 2^{j+1}\}}\neq0$$
if and only if at least one of the following two statements holds:

(i) $|y-z-h|\le 2^{j+1}$ and $|y-z|>2^{j+1}$;

(ii) $|y-z|\le 2^{j+1}$ and $|y-z-h|>2^{j+1}$.

\noindent This  together with the fact that $|h|\leq 2^j$ implies the following two cases:

(i) $2^{j+1}\leq|y-z|\leq2^{j+1}+|h|$;

(ii) $2^{j+1}-|h|\leq|y-z|\leq2^{j+1}$.

\noindent Since $|y-z|\ge 2|h|,$ then $|y-z|\simeq |z-y-h|$. Thus,
\begin{align*}
A_{5}
&\lesssim\frac{1}{2^{j(n+1)}}\int_{2^{j+1}-|h|\leq|y-z|\leq2^{j+1}+|h|}\frac{1}{|y-z|^{n-1}}\,dz\lesssim\frac{|h|}{2^{j(n+1)}}\chi_{\{|x-y|\lesssim2^{j}\}}.
\end{align*}
Combined with the estimates of $A_{1}$ and $A_2$, we finally obtain the desired estimate: for $|h|\le 2^j$ and some $\theta\in (0,1),$
\begin{align}\label{aj4}|a_j(x,y+h)-a_j(x,y)|\lesssim\bigg(\frac{|h|^\theta}{2^{j(n+\theta)}}
+\|\omega^{1/2}\|_{Dini}\frac{\omega^{1/2}(\frac{|h|}{2^j})}{2^{jn}}\bigg)\chi_{\{|x-y|\lesssim2^j\}}.\end{align}
Thus for $s\le 2^j,$ it follows from \eqref{aj4} by noting $Q_s1=0,$
\begin{align*}
|A_{j}Q_{s}f(x)|&=\bigg|\int_{{\Bbb R}^n}\int_{{\Bbb R}^n} (a_{j}(x,z)-a_{j}(x,y))\psi_{s}(z-y)dzf(y)\,dy\bigg|\\
&\lesssim \int_{|x-y|\lesssim2^j}\int_{|z-y|\le s}\frac{|z-y|^{\theta}}{2^{j(n+\theta)}}|\psi_{s}(z-y)|dz|f(y)|\,dy\\&\quad
+\|\omega^{1/2}\|_{Dini}\int_{|x-y|\lesssim2^j}\int_{{\Bbb R}^n}\frac{\omega^{1/2}(|z-y|/2^{j})}{2^{jn}}|\psi_{s}(z-y)|dz|f(y)|\,dy\\
&\lesssim \frac{s^{\theta}}{2^{j\theta}}Mf(x)
+\|\omega^{1/2}\|_{Dini}\omega^{1/2}(\frac{s}{2^{j}})Mf(x).
\end{align*}
Hence for $s\le 2^j,$ we also get for some $\theta\in (0,1)$
\begin{align}\label{tjh}
\|A_{j}Q_{s}f\|_{L^{2}}&\lesssim\bigg(\big(\frac{s}{2^{j}}\big)^{\theta}+\|\omega^{1/2}\|_{Dini}\omega^{1/2}(\frac{s}{2^{j}})\bigg)\|f\|_{L^{2}}.
\end{align}
 Together \eqref{A} with \eqref{tjh}, we verify for $A_j$ the assumption \eqref{zj} in Lemma \ref{gj}. Therefore
  we establish \eqref{I2}.\qed

\subsection{{Proof of $I_3$}}
For $I_3$, by the definition of the jump quantity and Chebychev's inequality,
\begin{align*}
\sup_{\lambda>0}\bigg\|\lambda\sqrt{N_\lambda(\{(\delta-\varphi_j)*T_jf\}_{j\in \Bbb Z})}\bigg\|_{L^2}\leq\bigg\|\bigg(\sum_{j\in\mathbb{Z}}|
(\delta-\varphi_j)*T_jf|^2\bigg)^{\frac{1}{2}}\bigg\|_{L^2}.
\end{align*} So to estimate $I_3$, it suffices to prove \begin{align} \label{I3}
\bigg\|\bigg(\sum_{j\in\mathbb{Z}}|
(\delta-\varphi_j)*T_jf|^2\bigg)^{\frac{1}{2}}\bigg\|_{L^2}\lesssim (1+\|\omega^{1/2}\|_{Dini}^2)\|f\|_{L^2}.
\end{align}
We will apply Lemma \ref{gj} to deal with the inequality \eqref{I3} by verifying the assumption  \eqref{zj} for $B_{j}=(\delta-\varphi_j)*T_j$.
Denote the kernel of $B_{j}$ as
$$b_{j}(x,y)=\int_{{\Bbb R}^n}\varphi_{j}(x-z)\big(K(x,y)\chi_{\{|x-y|>2^{j+1}\}}-K(z,y)\chi_{\{|z-y|>2^{j+1}\}}\big)\,dz.$$
Observe that $b_{j}$ is supported on the set $\{|x-y|>2^{j}\}$ and $\int_{{\Bbb R}^n}b_j(x,y)\,dy=0$ by $\int_{{\Bbb R}^n}K(\eta,y)\chi_{\{|\eta-y|>2^{j+1}\}}\,dy=\lim_{N\rightarrow \infty}\int_{|\eta-y|\le N}K(\eta,y)\chi_{\{|\eta-y|>2^{j+1}\}}\,dy=0.$
First, we deal with the case of $2^j\le s.$ Since $|x-z|\leq 2^j\leq\frac{|x-y|}{2},$  then by \eqref{k1} we get
\begin{align*}
|b_j(x,y)|
&\leq\int_{{\Bbb R}^n}|\varphi_{j}(x-z)||K(x,y)-K(z,y)|\chi_{\{|z-y|>2^{j+1}\}} dz\\&\quad +\int_{{\Bbb R}^n}|\varphi_{j}(x-z)||K(x,y)||\chi_{\{|x-y|>2^{j+1}\}}-\chi_{\{|z-y|>2^{j+1}\}}| dz
\\\nonumber&\lesssim \int_{{\Bbb R}^n}|\varphi_{j}(x-z)|\frac{\omega(|x-z|/|x-y|)}{|x-y|^{n}}dz\chi_{\{|x-y|> 2^{j}\}}+ \frac{1}{|x-y|^{n}}\chi_{\{2^{j}<|x-y|\leq 2^{j+2}\}}\int_{{\Bbb R}^n}|\varphi_{j}(x-z)|dz\\\nonumber
&\lesssim\bigg(\frac{\omega(2^j/|x-y|)}{|x-y|^{n}}+\frac{2^{j}}{|x-y|^{n+1}}\bigg)\chi_{\{|x-y|> 2^{j}\}}.
\end{align*}
where we have used that $|x-z|\le 2^j/2,$
 $ \chi_{\{|x-y|>2^{j+1}\}}-\chi_{\{|z-y|>2^{j+1}\}}$ is non-zero if
 $2^{j}<|x-y|\leq 2^{j+2}$. Write  $\omega_1(t)=\omega(t)+t$.
Again, by $\int_{{\Bbb R}^n} b_j(x,z)\,dz=0,$ we get
\begin{align*}
 |B_{j}Q_{s}f(x)|&=\bigg|\int_{{\Bbb R}^n}\int_{{\Bbb R}^n} b_{j}(x,z)\big(\psi_{s}(z-y)-\psi_{s}(x-y)\big)dzf(y)dy\bigg|\\
 &\lesssim\int_{|x-y|\lesssim s}\int_{2^j\leq |x-z|\leq s} \frac{\omega_1(2^j/|x-z|)}{|x-z|^{n}} |\psi_{s}(z-y)-\psi_{s}(x-y)| dz|f(y)|dy\\
 &\quad+\int_{{\Bbb R}^n}\int_{|x-z|>s}\frac{\omega_1(2^j/|x-z|)}{|x-z|^{n}}|\psi_{s}(z-y)|dz|f(y)|dy\\
  &\quad+\int_{{\Bbb R}^n}\int_{|x-z|>s}\frac{\omega_1(2^j/|x-z|)}{|x-z|^{n}}|\psi_{s}(x-y)|dz|f(y)|dy\\&=:B_1+B_2+B_3.
 \end{align*}
We will estimate $B_1,B_2$ and $B_3$ separately. First for $B_1,$ \begin{align*}
B_1&\lesssim\int_{|x-y|\lesssim s}\int_{2^j\leq |x-z|\leq s} \frac{\omega_1(2^j/|x-z|)}{|x-z|^{n}}\frac{|x-z|}{s^{n+1}}dz|f(y)|dy\\&\lesssim s^{-n}\int_{|x-y|\leq s}|f(y)|dy\frac{1}{s}\int_{2^j\le |x-z|\leq s} \frac{\omega_1(2^j/|x-z|)}{|x-z|^{n-1}} dz\\&\lesssim Mf(x)\omega_1^{1/2}(1)\frac{1}{s}\sum_{k=-\infty}^0\omega_1^{1/2}(2^{-k}2^j/s)\int_{ 2^{k-1}s<|x-z|\leq 2^ks} \frac{1}{|x-z|^{n-1}} dz.
\end{align*}
By the sub-additivity of $\omega_1(t)$, we get for any fixed $k\in \Bbb N_+,$  $\omega_1(kt)\le k\omega_1(t).$ Then we get
 \begin{align*}B_1&\lesssim \omega_1^{1/2}(1)\omega_1^{1/2}(2^j/s)Mf(x)\sum_{k=-\infty}^02^k2^{-k/2}\lesssim\|\omega_1^{1/2}\|_{Dini}\omega_1^{1/2}(2^j/s) Mf(x).
\end{align*}
For $B_2,$ we get
\begin{align*}
B_2&\lesssim\int_{|x-z|>s}\frac{\omega_1(2^j/|x-z|)}{|x-z|^{n}}Mf(z)dz\\
 &\lesssim\omega_1^{1/2}(\frac{2^j}{s})\sum_{k=0}^\infty\omega_1^{1/2}({2^{-k}})\int_{2^{k}s< |x-z|\le 2^{k+1}s}\frac{Mf(z)}{|x-z|^{n}}dz\\
 &\lesssim\|\omega_1^{1/2}\|_{Dini}\omega_1^{1/2}(\frac{2^j}{s})M^2f(x).
\end{align*}
For $B_3,$ similarly, we get
\begin{align*}
 B_3&\lesssim Mf(x)\int_{|x-z|>s}\frac{\omega_1(2^j/|x-z|)}{|x-z|^{n}}dz\lesssim\|\omega_1^{1/2}\|_{Dini}\omega_1^{1/2}(\frac{2^j}{s})Mf(x).
\end{align*}
Together the estimates of $B_1,\,B_2 $ and $B_3,$ we get
\begin{align}\label{sj0}
 \|B_{j}Q_{s}f\|_{L^2}&\lesssim \|\omega_1^{1/2}\|_{Dini}\omega_1^{1/2}(\frac{2^j}{s})\|f\|_{L^2}\lesssim (\omega^{1/2}(s/2^j)\|\omega^{1/2}\|_{Dini}+(\frac{s}{2^{j}})^{1/2})\|f\|_{L^2}.
\end{align}
Secondly, we deal with the case of $s\le 2^j.$
For $|h|\le 2^j,$
\begin{align*}
&|b_{j}(x,y+h)-b_{j}(x,y)|\\
&\leq\bigg|\int_{{\Bbb R}^n}\varphi_{j}(x-z)\big(K_{j}(x,y+h)-K_{j}(x,y)\big)\,dz\bigg|+\bigg|\int_{{\Bbb R}^n}\varphi_{j}(x-z)\big(K_{j}(z,y+h)-K_{j}(z,y)\big)\,dz\bigg|\\
&=:B_4+B_5.
\end{align*}
For $B_4$, we get
\begin{align*}
B_4
&\lesssim\int_{{\Bbb R}^n}|\varphi_{j}(x-z)||K(x,y+h)-K(x,y)|\chi_{\{|x-y-h|>2^{j+1}\}}\,dz\\
&\quad+\int_{{\Bbb R}^n}|\varphi_{j}(x-z)||K(x,y)||\chi_{\{|x-y-h|>2^{j+1}\}}-\chi_{\{|x-y|>2^{j}\}}|\,dz\\
&\lesssim \frac{\omega(|h|/|x-y|)}{|x-y|^{n}}\chi_{\{|x-y|>2^{j}\}}+\frac{1}{|x-y|^n}\int_{{\Bbb R}^n}|\varphi_{j}(x-z)||\chi_{\{|x-y-h|>2^{j+1}\}}-\chi_{\{|x-y|>2^{j+1}\}}|\,dz,
\end{align*}
where we used the fact $|h|\le2^j$ and $|x-y-h|\ge2^{j+1}$ imply $|x-y|>2^j$.
For the second quantity of the right hand side, by $|h|\le2^j$, then $\chi_{\{|x-y-h|>2^{j+1}\}}-\chi_{\{|x-y|>2^{j+1}\}}$ is nonzero  if
 $2^{j+1}-|h|\leq|x-y|\le 2^{j+1}+|h|.$
Therefore, one has
\begin{align*}B_4\lesssim \frac{\omega(|h|/|x-y|)}{|x-y|^{n}}\chi_{\{|x-y|>2^{j}\}}+\frac{1}{|x-y|^n}\chi_{\{2^{j+1}-|h|\leq|x-y|\le2^{j+1}+|h|\}}.\end{align*}
Similarly, for $B_5$,
\begin{align*}
B_5
&\lesssim\int_{{\Bbb R}^n}|\varphi_{j}(x-z)|\frac{\omega(|h|/|z-y|)}{|z-y|^{n}}\chi_{\{|z-y|>2^{j}\}}\,dz+\int_{{\Bbb R}^n}\frac{|\varphi_{j}(x-z)|}{|z-y|^n}\chi_{\{2^{j+1}-|h|\le|z-y|\leq2^{j+1}+|h|\}}\,dz.
\end{align*}
For $|x-z|\le2^j/2$ and $|z-y|>2^j$, then
$
2|z-y|\ge |x-y|\geq\frac{|z-y|}{2}>2^j/2.
$
We get for some $\theta\in (0,1)$ and $\frac{1}{p'}=\theta,$
\begin{align*}
B_5&\lesssim \frac{\omega(|h|/|x-y|)}{|x-y|^{n}}\chi_{\{|x-y|\ge2^{j}/2\}}+\frac{1}{|x-y|^{n+\theta}}\chi_{\{|x-y|\ge2^{j}/2\}} \bigg(\int_{{\Bbb R}^n}|\varphi_j(x-z)|^p\,dz\bigg)^{\frac{1}{p}}\\&\quad\times\bigg(\int_{2^{j+1}-|h|\leq|z-y|\le2^{j+1}+|h|}|z-y|^{\theta p'}\,dz\bigg)^{\frac{1}{p'}}\\
&\lesssim\frac{\omega(|h|/|x-y|)}{|x-y|^{n}}\chi_{\{|x-y|\ge2^{j}/2\}}+\frac{|h|^\theta}{|x-y|^{n+\theta}}\chi_{\{|x-y|\ge2^{j}/2\}}.
\end{align*}
Combined the estimate of $B_4$ and $B_5$, we get for $|h|\le 2^j$
\begin{align}\label{bjy}
|b_j(x,y+h)-b_j(x,y)|&\lesssim\bigg(\frac{\omega(|h|/|x-y|)}{|x-y|^{n}}+
\frac{|h|^\theta}{|x-y|^{n+\theta}}\bigg)\chi_{\{|x-y|\ge2^{j}/2\}}+\frac{1}{|x-y|^n}\chi_{\{2^{j+1}-|h|\leq|x-y|\le2^{j+1}+|h|\}}
\end{align}
 Now we apply \eqref{bjy} to estimate  $\|B_jQ_sf\|_{L^2}.$ Since $Q_s1=0,$ and $s\leq 2^{j}$,
\begin{align*}
|B_{j}Q_{s}f(x)|
&=\bigg|\int_{{\Bbb R}^n}\int_{{\Bbb R}^n} \big(b_{j}(x,z)-b_{j}(x,y)\big)\psi_{s}(z-y)dzf(y)dy\bigg|\\
&\lesssim\int_{{\Bbb R}^n}\int_{|x-y|\ge 2^j/2\bigcap|z-y|\le s} \frac{\omega(|z-y|/|x-y|)}{|x-y|^{n}}|\psi_{s}(z-y)|dz|f(y)|dy\\
&\quad+\int_{{\Bbb R}^n}\int_{\{2^{j+1}-|z-y|\leq|x-y|\leq2^{j+1}+|z-y|\}\bigcap\{|z-y|\le s\}} |\psi_{s}(z-y)|dz\frac{|f(y)|}{|x-y|^{n}}dy\\&\lesssim\int_{|x-y|\ge 2^j/2} \frac{2\omega(\frac{s}{2|x-y|})}{|x-y|^{n}}|f(y)|dy+\bigg(\int_{2^{j+1}-s\leq|x-y|\leq2^{j+1}+s} \frac{1}{|x-y|^{p'n}}dy\bigg)^{1/p'}\bigg(\int_{|x-y|\leq2^{j+2}} |f(y)|^{p}dy\bigg)^{1/p}\\&\lesssim\omega^{1/2}(s/2^j)\int_{|x-y|\ge 2^j/2} \frac{\omega^{1/2}(\frac{s}{2|x-y|})}{|x-y|^{n}}|f(y)|dy+(\frac{s}{2^{j}})^{\frac{1}{p'}}\bigg(M|f|^{p}(x)\bigg)^{1/p}
\\&\lesssim\omega^{1/2}(s/2^j)\|\omega^{1/2}\|_{Dini}Mf(x)+(\frac{s}{2^{j}})^{\frac{1}{p'}}\bigg(M|f|^{p}(x)\bigg)^{1/p}.
\end{align*}
where  $1<p<2.$
Therefore we get for $s\leq 2^j$ and some  $\theta\in (0,1)$
\begin{align}\label{sj1}
\|B_{j}Q_{s}f\|_{L^{2}}\lesssim(\omega^{1/2}(s/2^j)\|\omega^{1/2}\|_{Dini}+(\frac{s}{2^{j}})^\theta)\|f\|_{L^{2}}.
\end{align}
Therefore, combined  with \eqref{sj0}, we verify the assumption  \eqref{zj} in Lemma \ref{gj}.
 Therefore we establish \eqref{I3}.

\section{Proof of Theorem \ref{tJ2}---The short variational estimate \eqref{short variation singular}} \label{s6}
 Let $T_{j,\,r}f(x)=\int_{2^{j}r\le |x|\leq2^{j+1}} K(x,y)f(y)\,dy,$ $r\in [1,2]$. Then
\begin{align*}
S_{2}(\mathcal{T} f)(x)
&=\Big(\dsum_{j\in\mathbb{Z}}\dsup_{\substack
{t_1<\cdots<t_N\\
[t_l,t_{l+1}]\subset[1,2]}}\dsum_{l=1}^{N-1}|T_{j,\,t_l}
f(x)-T_{j,\,t_{l+1}}
f(x)|^2\Big)^{\frac{1}{2}}\\
&= \Big(\dsum_{j\in\mathbb{Z}}\dsup_{\substack
{t_1<\cdots<t_N\\
[t_l,\,t_{l+1}]\subset[1,2]}}\dsum_{l=1}^{N-1}|{T}_{j,\,t_l,\,t_{l+1}}
f(x)|^2\Big)^{\frac{1}{2}},\end{align*}
where the operator ${{T}}_{j,\,t_l,\,t_{l+1}}$ is given by
 $${{T}}_{j,\,t_l,\,t_{l+1}}f(x)=\dint_{2^jt_l<|x-y|\le 2^{j}t_{l+1}}K(x,y)f(y)dy,\,\,[t_l,t_{l+1}]\subset[1,2].$$
Then by the Calder\'on identity $\int_0^\infty Q_s^2\frac{ds}{s}=\mathcal{I},$ one has
\begin{align*}
\|S_{2}(\mathcal{T}f)\|_{L^2}^2
&=\sum_{j\in\mathbb{Z}}\int_{{\Bbb R}^n}\dsup_{\substack
{t_1<\cdots<t_N\\
[t_l,t_{l+1}]\subset[1,2]}}\dsum_{l=1}^{N-1}\bigg|T_{j,\,t_l,\,t_{l+1}}\int_0^\infty Q_s^2f(x)\frac{ds}{s}\bigg|^2\,dx\\&\lesssim \sum_{j\in\mathbb{Z}}\int_{{\Bbb R}^n}\dsup_{\substack
{t_1<\cdots<t_N\\
[t_l,t_{l+1}]\subset[1,2]}}\dsum_{l=1}^{N-1}\bigg|\int_{2^{j-1}}^\infty T_{j,\,t_l,\,t_{l+1}}Q_s Q_sf(x)\frac{ds}{s}\bigg|^2\,dx\\&\quad+\sum_{j\in\mathbb{Z}}\int_{{\Bbb R}^n}\dsup_{\substack
{t_1<\cdots<t_N\\
[t_l,t_{l+1}]\subset[1,2]}}\dsum_{l=1}^{N-1}\bigg|\int_0^{2^{j-1}} T_{j,\,t_l,\,t_{l+1}} Q_s Q_sf(x)\frac{ds}{s}\bigg|^2\,dx\\&=: I+II.
\end{align*}
For $I,$
denote by $K_{j,\, l,\, s}(x,y)$ the kernel of $T_{j,\, t_l,\, t_{l+1}}Q_s.$ Since $T_{j,\, t_l,\, t_{l+1}}1=0$ for any $j\in \Bbb Z,$ then
\begin{align*}
\dsum_{l=1}^{N-1} |K_{j,\, l,\, s}(x,y)|&=\dsum_{l=1}^{N-1}\bigg|\int_{2^jt_l<|x-y|\le 2^{j}t_{l+1}} K(x,z)\big(\psi_{s}(z-y)-\psi_{s}(x-y)\big)dz\bigg|\\&\le \int_{2^jt_1<|x-z|\le 2^{j}t_N} |K(x,z)||\psi_{s}(z-y)-\psi_{s}(x-y)|dz.
\end{align*}
First,  by \eqref{k0} and $t_1,\,t_N\in [1,2]$, we get $|K(x,z)|\chi_{\{2^jt_1<|x-z|\le 2^{j}t_N\}}\lesssim \frac{2^{j/2}}{|x-z|^{n+1/2}},$ then
\begin{align*}
\dsum_{l=1}^{N-1}|T_{j,\, t_l,\, t_{l+1}} Q_sf(x)|&\le\int_{|x-y|\leq4s}\int_{|x-z|\le s} \frac{2^{j/2}}{|x-z|^{n+1/2}}|\psi_{s}(z-y)-\psi_{s}(x-y)|dz|f(y)|dy\\&\quad+\int_{{\Bbb R}^n}\int_{|x-z|> s} \frac{2^{j/2}}{|x-z|^{n+1/2}}|\psi_{s}(z-y)|dz|f(y)|dy\\&\quad+\int_{{\Bbb R}^n}\int_{|x-z|> s} \frac{2^{j/2}}{|x-z|^{n+1/2}}|\psi_{s}(x-y)|dz|f(y)|dy\\
&\lesssim\int_{|x-y|\leq4s}\int_{|x-z|\le s}\frac{2^{j/2}}{|x-z|^{n+1/2}}\frac{|x-z|}{s^{n+1}}dz|f(y)|dy+(\frac{2^{j}}{s})^{1/2}(M^2f(x)+Mf(x))\\
&\lesssim(\frac{2^{j}}{s})^{1/2}(Mf(x)+M^2f(x)).
\end{align*}

Then by the H\"{o}lder's inequality
\begin{align*}
I
&\lesssim \sum_{j\in\mathbb{Z}}\int_{{\Bbb R}^n}\bigg(\int_{2^{j-1}}^\infty (\frac{2^{j}}{s})^{1/2}(MQ_sf(x)+M^2Q_sf(x)) \frac{ds}{s}\bigg)^2\,dx\\&\lesssim \sum_{j\in\mathbb{Z}}\int_{{\Bbb R}^n}\bigg(\int_{2^{j-1}}^\infty (\frac{2^{j}}{s})^{1/2} \frac{ds}{s}\bigg)\bigg(\int_{2^{j-1}}^\infty (\frac{2^{j}}{s})^{1/2}\big(MQ_sf(x)+M^2Q_sf(x)\big)^2 \frac{ds}{s}\bigg)\,dx\\&\lesssim \int_{0}^\infty\int_{{\Bbb R}^n} \big(\sum_{2^{j-1}\le s}(\frac{2^{j}}{s})^{1/2}\big)(MQ_sf(x)+M^2Q_sf(x))^2\,dx \frac{ds}{s}\\&\lesssim \int_{{\Bbb R}^n}\int_{0}^\infty |Q_sf(x)|^2 \frac{ds}{s}\,dx\\&\lesssim \|f\|_{L^2}^2.
\end{align*}
Next, for $II$,
since $Q_{s}1=0$,  then
\begin{align*}
\dsum_{l=1}^{N-1}|T_{j,\,t_l,\,t_{l+1}}Q_{s}f(x)|^2
&=\dsum_{l=1}^{N-1}\bigg|\int_{{\Bbb R}^n}\int_{{\Bbb R}^n}\big(K(x,z)\chi_{\{2^jt_l<|x-z|\le 2^{j}t_{l+1}\}}-K(x,y)\chi_{\{2^jt_l<|x-y|\le 2^{j}t_{l+1}\}}\big)\psi_{s}(z-y)\,dzf(y)\,dy\bigg|^2\\&\le\bigg(\dsum_{l=1}^{N-1}\int_{{\Bbb R}^n}\int_{2^jt_l< |x-z|\le 2^jt_{l+1}}|K(x,z)-K(x,y)||\psi_{s}(z-y)|\,dz|f(y)|\,dy\bigg)^2\\&\quad+\dsum_{l=1}^{N-1}\bigg(\int_{{\Bbb R}^n}\int_{{\Bbb R}^n}|K(x,y)||\chi_{\{2^jt_l< |x-z|\le 2^jt_{l+1}\}}-\chi_{\{2^jt_l< |x-y|\le 2^jt_{l+1}\}}||\psi_{s}(z-y)|\,dz|f(y)|\,dy\bigg)^2\\&=:II_1+II_2,
\end{align*}
 For $II_1, $ by $s\le 2^{j-1}$  and \eqref{k1} we get
 \begin{align*}
II_1&\le\bigg(\int_{{\Bbb R}^n}\int_{2^jt_1< |x-z|\le 2^jt_{N}}|K(x,z)-K(x,y)||\psi_{s}(z-y)|\,dz|f(y)|\,dy\bigg)^2\\&\lesssim\bigg(\int_{{\Bbb R}^n}\int_{{\Bbb R}^n}\frac{\omega(|z-y|/|x-y|)}{|x-y|^n}\chi_{\{2^{j}/2\le |x-y|\le 2^{j+2}\}}|\psi_{s}(z-y)|\,dz|f(y)|\,dy\bigg)^2\\
&\lesssim \omega^2 (\frac{s}{2^j})(Mf(x))^2.
\end{align*}
 For $II_2.$ Note that
\begin{align*}
\chi_{\{2^{j}t_l<|x-z|\leq2^{j}t_{l+1}\}}-\chi_{\{2^{j}t_l<|x-y|\leq2^{j}t_{l+1}\}}\neq0,
\end{align*}
if and only if at least one of the following four statements holds,

(i)\, $2^{j}t_l<|x-z|\leq2^{j}t_{l+1}$ and $|x-y|\le2^{j}t_l$;

(ii)\, $2^{j}t_l<|x-z|\leq2^{j}t_{l+1}$ and $|x-y|>2^{j}t_{l+1}$;

(iii)\, $2^{j}t_l<|x-y|\leq2^{j}t_{l+1}$ and $|x-z|\le2^{j}t_l$;

(iv)\, $2^{j}t_l<|x-y|\leq2^{j}t_{l+1}$ and $|x-z|>2^{j}t_{l+1}$.

\noindent This together with the fact that $|y-z|\leq s$ implies the following four cases

(i)\, $2^{j}t_l<|x-z|\leq2^{j}t_{l+1}$ and $2^{j}t_l<|x-z|\leq2^{j}t_l+s;$

(ii)\, $2^{j}t_l<|x-z|\leq2^{j}t_{l+1}$and $2^{j}t_{l+1}-s\leq|x-z|\leq2^{j}t_{l+1};$

(iii) \, $2^{j}t_l<|x-y|\leq2^{j}t_{l+1}$and $2^{j}t_l<|x-y|\leq2^{j}t_l+s;$

(iv)\, $2^{j}t_l<|x-y|\leq2^{j}t_{l+1}$ and $2^{j}t_{l+1}-s\leq|x-y|\leq2^{j}t_{l+1}.$

\noindent Case (i) and Case (ii )  can be dealt with similarly. We only consider Case (i).   Taking an arbitrary $1<p<2,$ by \eqref{k0} and the H\"{o}lder inequality, we get \begin{align*}
II_2
&\lesssim\dsum_{l=1}^{N-1}\bigg(\frac{1}{2^{jn}}\int_{2^{j}t_l<|x-z|\leq2^{j}t_{l+1}\bigcap 2^{j}t_l<|x-z|\leq2^{j}t_l+s}
\int_{{\Bbb R}^n}|\psi_{s}(z-y)||f(y)|\,dy\,dz\bigg)^{2}\\
&\lesssim \dsum_{l=1}^{N-1}\bigg(\frac{1}{2^{jn}}\int_{2^{j}t_l<|x-z|\leq2^{j}t_{l+1}\bigcap 2^{j}t_l<|x-z|\leq2^{j}t_l+s}
Mf(z)\,dz\bigg)^{2}\\
&\lesssim \dsum_{l=1}^{N-1}\bigg(\frac{1}{2^{jn}}\int_{2^{j}t_l<|x-z|\leq2^{j}t_{l+1}}
Mf(z)^p\,dz\bigg)^{2/p}\bigg(\frac{1}{2^{jn}}\int_{2^{j}t_l<|x-z|\leq2^{j}t_l+s}
\,dz\bigg)^{2/p'}\\
&\lesssim \bigg(\frac{1}{2^{jn}}\dsum_{l=1}^{N-1}\int_{2^{j}t_l<|x-z|\leq2^{j}t_{l+1}}
Mf(z)^p\,dz\bigg)^{2/p}(\frac{s}{2^j})^{2/p'}\\
&\lesssim \bigg(\frac{1}{2^{jn}}\int_{2^{j}t_1<|x-z|\leq2^{j}t_{N}}
Mf(z)^p\,dz\bigg)^{2/p}(\frac{s}{2^j})^{2/p'}\\
&\lesssim
(M(Mf)^p)^{2/p}(x)(\frac{s}{2^j})^{2/p'}.
\end{align*}

\noindent Case (iii) and Case (iv) can be treated similarly. We only consider Case (iii).   Taking $1<p<2,$ by \eqref{k0} and the H\"{o}lder inequality, we get \begin{align*}
II_2
&\lesssim\dsum_{l=1}^{N-1}\bigg(\frac{1}{2^{jn}}\int_{\{2^{j}t_l<|x-y|\leq2^{j}t_{l+1}\}\bigcap\{ 2^{j}t_l<|x-y|\leq2^{j}t_{l}+s\}}|f(y)|\,dy\bigg)^2\\
&\lesssim \dsum_{l=1}^{N-1}\bigg(\frac{1}{2^{jn}}\int_{2^{j}t_l<|x-y|\leq2^{j}t_{l+1}}
|f(y)|^p\,dy\bigg)^{2/p}\bigg(\frac{1}{2^{jn}}\int_{2^{j}t_l<|x-y|\leq2^{j}t_l+s}
\,dy\bigg)^{2/p'}\\
&\lesssim \bigg(\frac{1}{2^{jn}}\dsum_{l=1}^{N-1}\int_{2^{j}t_l<|x-y|\leq2^{j}t_{l+1}}
|f(y)|^p\,dy\bigg)^{2/p}(\frac{s}{2^j})^{2/p'}\\
&\lesssim \bigg(\frac{1}{2^{jn}}\int_{2^{j}t_1\leq|x-y|\leq2^{j}t_{N}}
|f(y)|^p\,dy\bigg)^{2/p}(\frac{s}{2^j})^{2/p'}\\
&\lesssim
(M|f|^p)^{2/p}(x)(\frac{s}{2^j})^{2/p'}.
\end{align*}
To conclude the estimate of $II_1$  and $II_2$ we get for some $\theta\in (0,1),\,1<p<2,
$\begin{align}\label{II}
\dsum_{l=1}^{N-1}|T_{j,t_l,t_{l+1}}Q_{s}f(x)|^2&\lesssim\big(\omega(\frac{s}{2^j})+\big(\frac{s}{2^{j}}\big)^{\theta}\big)^2(M(Mf)^p)^{2/p}(x)\simeq \widetilde{\omega}^2(\frac{s}{2^j})(M(Mf)^p)^{2/p}(x).
\end{align}
Here $\widetilde{\omega}(t)=\omega(t)+t^\theta.$ Now by \eqref{II} we can conclude the estimate of $II$,
 \begin{align*}
II
&\le \sum_{j\in\mathbb{Z}}\int_{{\Bbb R}^n}\dsup_{\substack
{t_1<\cdots<t_N\\
[t_l,t_{l+1}]\subset[1,2]}}\dsum_{l=1}^{N-1}\bigg|\int_{0}^{2^{j-1}}\widetilde{\omega}^{1/2}(\frac{s}{2^j}) \widetilde{\omega}^{-1/2}(\frac{s}{2^j})|T_{j,t_l,t_{l+1}}Q_s Q_sf(x)|\frac{ds}{s}\bigg|^2\,dx\\&\le \sum_{j\in\mathbb{Z}}\int_{{\Bbb R}^n}\dsup_{\substack
{t_1<\cdots<t_N\\
[t_l,t_{l+1}]\subset[1,2]}}\bigg(\int_{0}^{2^{j-1}}\widetilde{\omega}(\frac{s}{2^j})\frac{ds}{s} \bigg)\bigg(\int_{0}^{2^{j-1}}\widetilde{\omega}^{-1}(\frac{s}{2^j})\dsum_{l=1}^{N-1}|T_{j,t_l,t_{l+1}}Q_s Q_sf(x)|^2\frac{ds}{s}\bigg)\,dx\\&\lesssim \|\widetilde{\omega}\|_{Dini}\int_{0}^{\infty}\int_{{\Bbb R}^n}\big(\sum_{2^{j-1}\ge s}\widetilde{\omega}(\frac{s}{2^j})\big)(M(MQ_sf)^p)^{2/p}(x)\,dx\frac{ds}{s}\\&\lesssim \|\widetilde{\omega}\|_{Dini}^2\int_{{\Bbb R}^n}\int_{0}^{\infty}|Q_sf(x)|^2\frac{ds}{s}\,dx\\&\lesssim(\|\omega\|_{Dini}+1)^2\|f\|_{L^2}^2.
\end{align*}

 Combining the estimates of $I$ and $II,$  we get
\begin{align*}
\|S_{2}(\mathcal{T}f)\|_{L^2}^2\lesssim(1+\|\omega\|_{Dini})^2\|f\|_{L^2}^2.
\end{align*}
\qed

\end{document}